\expandafter\ifx\csname mthreemacsloaded\endcsname\relax\else \fi

\magnification1100
\input amstex


 \catcode`\@=11
 \let\wlog@ld\wlog
 \def\wlog#1{\relax}

 \newif\ifIN@
 \def\m@rker{\m@@rker}
 \def\IN@{\expandafter\INN@\expandafter}
 \long\def\INN@0#1@#2@{\long\def\NI@##1#1##2##3\ENDNI@
    {\ifx\m@rker##2\IN@false\else\IN@true\fi}%
     \expandafter\NI@#2@@#1\m@rker\ENDNI@}
  \newtoks\Initialtoks@  \newtoks\Terminaltoks@
  \def\SPLIT@{\expandafter\SPLITT@\expandafter}
  \def\SPLITT@0#1@#2@{\def\TTILPS@##1#1##2@{%
     \Initialtoks@{##1}\Terminaltoks@{##2}}\expandafter\TTILPS@#2@}
  \newtoks\Trimtoks@

 \def\ForeTrim@{\expandafter\ForeTrim@@\expandafter}
 \def\ForePrim@0 #1@{\Trimtoks@{#1}}
 \def\ForeTrim@@0#1@{\IN@0\m@rker. @\m@rker.#1@%
     \ifIN@\ForePrim@0#1@%
     \else\Trimtoks@\expandafter{#1}\fi}
 
  \def\Trim@0#1@{%
      \ForeTrim@0#1@%
      \IN@0 @\the\Trimtoks@ @%
        \ifIN@
             \SPLIT@0 @\the\Trimtoks@ @\Trimtoks@\Initialtoks@
             \IN@0\the\Terminaltoks@ @ @%
                 \ifIN@
                 \else \Trimtoks@ {FigNameWithSpace}%
                 \fi
        \fi
      }

  \font\titlebold=cmbx12 scaled 1200
  \font\twelvebold=cmbx12
  \font\tenbold=cmbx10
  \font\ninebold=cmbx9
  \font\sevenbold=cmbx7
  \font\fivebold=cmbx5

  \input amssym.def \input amssym
     \font\titlemsa=msam10 at 14.4pt
     \font\titlemsb=msbm10 at 14.4pt
     \font\titleeufm=eufm10 at 14.4pt
     \font\twelvemsa=msam10 scaled 1200
     \font\twelvemsb=msbm10 scaled 1200
     \font\twelveeufm=eufm10 scaled 1200
     \font\ninemsa=msam9
     \font\ninemsb=msbm9
     \font\nineeufm=eufm9

   \ifx\cyrfam\undefined
   \else
     \immediate\write16{}%
     \message{ !!! cyr fonts already defined. !!! }
     \message{ --- edit out superfluous font defs? }
   \fi
   \newfam\cyrfam
       \font\titlecyr=wncyr10 scaled 1440 
       \font\twelvecyr=wncyr10 scaled 1200
       \font\tencyr=wncyr10
       \font\ninecyr=wncyr9
       \font\sevencyr=wncyr7
       \font\sixcyr=wncyr6

   \newfam\eusmfam
       \font\titleeusm=eusm10 scaled 1440
       \font\twelveeusm=eusm10 scaled 1200
       \font\teneusm=eusm10
       \font\nineeusm=eusm9
       \font\seveneusm=eusm7
       
       \font\fiveeusm=eusm5

\let\Cal\cal

    \font\ninemrm=cmr9 
    \font\ninei=cmmi9
    \font\ninesy=cmsy9 
    \skewchar\ninei='177
    \skewchar\ninesy='60

  \font\twelvemrm=cmr10 at 12pt 
  \font\twelvei=cmmi10 at 12pt
  \font\twelvesy=cmsy10 at 12pt

  \font\titlemrm=cmr10 at 14.4pt 
  \font\titlei=cmmi10 at 14.4pt
  \font\titlesy=cmsy10 at 14.4pt


  \def\Smallfonts{\ninepoint}

  \def\Hfont{\titlepoint\bf}
  \def\Authorfont{\twelvepoint\it}
  \def\HHfont{\twelvepoint\bf}
  \def\HHHfont{\bf}
  \def\Bibfont{\tenbf}
  \def\Coordfont{\nineit }

  \def \thfont {\bf }
  \def \pffont {\it\itSpacing }
  \def \rkfont {\bf }
  \def \dffont {\bf }
  \def \egfont {\bf }

 \def\ninepoint{%
  \def\rm{\fam0\ninerm}%
    \textfont0=\ninemrm  \scriptfont0=\sevenrm  \scriptscriptfont0=\fiverm
    \textfont1=\ninei    \scriptfont1=\seveni   \scriptscriptfont1=\fivei
  \def\mit{\fam1\ninei}%
  \def\oldstyle{\fam1\ninei}%
    \textfont2=\ninesy   \scriptfont2=\sevensy  \scriptscriptfont2=\fivesy
    \textfont3=\tenex    \scriptfont3=\tenex    \scriptscriptfont3=\tenex
  \def\it{\fam\itfam\nineit}%
    \textfont\itfam=\nineit
  \def\bf{\ifmmode\fam\bffam\else\ninebf\fi}%
    \textfont\bffam=\ninebold 
    \scriptfont\bffam=\sevenbold 
    \scriptscriptfont\bffam=\fivebold%
  \def\msa{\fam\msafam\ninemsa}%
    \textfont\msafam=\ninemsa 
    \scriptfont\msafam=\sevenmsa
    \scriptscriptfont\msafam=\fivemsa%
  \def\msb{\fam\msbfam\ninemsb}%
    \textfont\msbfam=\ninemsb%
    \scriptfont\msbfam=\sevenmsb%
    \scriptscriptfont\msbfam=\fivemsb%
  \def\eufm{\fam\eufmfam\nineeufm}%
    \textfont\eufmfam=\nineeufm
    \scriptfont\eufmfam=\seveneufm
    \scriptscriptfont\eufmfam=\fiveeufm
   \def\eusm{\fam\eusmfam\nineeusm}%
     \textfont\eusmfam=\nineeusm
     \scriptfont\eusmfam=\seveneusm
     \scriptscriptfont\eusmfam=\fiveeusm
   \def\cyr{\fam\cyrfam\ninecyr}%
     \textfont\cyrfam=\ninecyr
     \scriptfont\cyrfam=\sevencyr
     \scriptscriptfont\cyrfam=\sixcyr
  \setbox\strutbox=\hbox{\vrule
      height7pt depth3pt width0pt}%
   \baselineskip=10.8pt\rm}

 \let\eightpoint\ninepoint 

 \def\tenpoint{%
  \def\rm{\fam0\tenrm}%
    \textfont0=\tenmrm \scriptfont0=\sevenrm \scriptscriptfont0=\fiverm%
  \def\mit{\fam1\teni}%
  \def\oldstyle{\fam1\teni}%
    \textfont1=\teni   \scriptfont1=\seveni  \scriptscriptfont1=\fivei%
    \textfont2=\tensy  \scriptfont2=\sevensy \scriptscriptfont2=\fivesy%
    \textfont3=\tenex  \scriptfont3=\tenex   \scriptscriptfont3=\tenex%
  \def\it{\fam\itfam\tenit}%
    \textfont\itfam=\tenit%
  \def\bf{\ifmmode\fam\bffam\else\tenbf\fi}%
    \textfont\bffam=\tenbold
    \scriptfont\bffam=\sevenbold%
    \scriptscriptfont\bffam=\fivebold%
  \def\msa{\fam\msafam\tenmsa}%
    \textfont\msafam=\tenmsa%
    \scriptfont\msafam=\sevenmsa%
    \scriptscriptfont\msafam=\fivemsa%
  \def\msb{\fam\msbfam\tenmsb}%
    \textfont\msbfam=\tenmsb%
    \scriptfont\msbfam=\sevenmsb%
    \scriptscriptfont\msbfam=\fivemsb%
  \def\eufm{\fam\eufmfam\teneufm}%
   \textfont\eufmfam=\teneufm
   \scriptfont\eufmfam=\seveneufm
   \scriptscriptfont\eufmfam=\fiveeufm
   \def\eusm{\fam\eusmfam\teneusm}%
    \textfont\eusmfam=\teneusm
    \scriptfont\eusmfam=\seveneusm
    \scriptscriptfont\eusmfam=\fiveeusm
   \def\cyr{\fam\cyrfam\tencyr}%
    \textfont\cyrfam=\tencyr
    \scriptfont\cyrfam=\sevencyr
    \scriptscriptfont\cyrfam=\sixcyr
  \setbox\strutbox=\hbox{\vrule %
      height8.5pt depth3.5ptwidth0pt}%
  \baselineskip=\StdBaselineskip\rm}

 \def\twelvepoint{%
  \def\rm{\fam0\twelverm}%
    \textfont0=\twelvemrm \scriptfont0=\tenmrm \scriptscriptfont0=\sevenrm
    \textfont1=\twelvei   \scriptfont1=\teni   \scriptscriptfont1=\seveni
  \def\mit{\fam1\twelvei}%
  \def\oldstyle{\fam1\twelvei}%
    \textfont2=\twelvesy  \scriptfont2=\tensy  \scriptscriptfont2=\sevensy
    \textfont3=\tenex  \scriptfont3=\tenex  \scriptscriptfont3=\tenex
  \def\it{\fam\itfam\twelveit}%
    \textfont\itfam=\twelveit
  \def\bf{\ifmmode\fam\bffam\else\twelvebf\fi}%
    \textfont\bffam=\twelvebold
    \scriptfont\bffam=\tenbold%
    \scriptscriptfont\bffam=\sevenbold%
  \def\msa{\fam\msafam\twelvemsa}%
    \textfont\msafam=\twelvemsa%
    \scriptfont\msafam=\tenmsa%
    \scriptscriptfont\msafam=\sevenmsa%
  \def\msb{\fam\msbfam\twelvemsb}%
    \textfont\msbfam=\twelvemsb%
    \scriptfont\msbfam=\tenmsb%
    \scriptscriptfont\msbfam=\sevenmsb%
  \def\eufm{\fam\eufmfam\twelveeufm}%
   \textfont\eufmfam=\twelveeufm
   \scriptfont\eufmfam=\teneufm
   \scriptscriptfont\eufmfam=\seveneufm
   \def\eusm{\fam\eusmfam\twelveeusm}%
    \textfont\eusmfam=\twelveeusm
    \scriptfont\eusmfam=\teneusm
    \scriptscriptfont\eusmfam=\seveneusm
   \def\cyr{\fam\cyrfam\tencyr}%
    \textfont\cyrfam=\twelvecyr
    \scriptfont\cyrfam=\tencyr
    \scriptscriptfont\cyrfam=\sevencyr
  \setbox\strutbox=\hbox{\vrule
      height10.2pt depth4.55pt width0pt}%
  \baselineskip=14pt\rm}

 \def\titlepoint{%
    \textfont0=\titlemrm \scriptfont0=\twelvemrm \scriptscriptfont0=\tenmrm
    \textfont1=\titlei   \scriptfont1=\twelvei   \scriptscriptfont1=\teni
  \def\mit{\fam1\titlei}%
  \def\oldstyle{\fam1\titlei}%
    \textfont2=\titlesy  \scriptfont2=\twelvesy  \scriptscriptfont2=\tensy
    \textfont3=\tenex
    \scriptfont3=\tenex
    \scriptscriptfont3=\tenex
  \def\it{\fam\itfam\titleit}%
    \textfont\itfam=\titleit
  \def\bf{\ifmmode\fam\bffam\else\titlebf\fi}%
    \textfont\bffam=\titlebold
    \scriptfont\bffam=\twelvebold%
    \scriptscriptfont\bffam=\tenbold%
  \def\msa{\fam\msafam\titlemsa}%
    \textfont\msafam=\titlemsa%
    \scriptfont\msafam=\twelvemsa%
    \scriptscriptfont\msafam=\tenmsa%
  \def\msb{\fam\msbfam\titlemsb}%
    \textfont\msbfam=\titlemsb%
    \scriptfont\msbfam=\twelvemsb%
    \scriptscriptfont\msbfam=\tenmsb%
  \def\eufm{\fam\eufmfam\titleeufm}%
    \textfont\eufmfam=\titleeufm
    \scriptfont\eufmfam=\twelveeufm
    \scriptscriptfont\eufmfam=\teneufm
   \def\eusm{\fam\eusmfam\titleeusm}%
     \textfont\eusmfam=\titleeusm
     \scriptfont\eusmfam=\twelveeusm
     \scriptscriptfont\eusmfam=\teneusm
   \def\cyr{\fam\cyrfam\tencyr}%
    \textfont\cyrfam=\titlecyr
    \scriptfont\cyrfam=\twelvecyr
    \scriptscriptfont\cyrfam=\tencyr
  \setbox\strutbox=\hbox{\vrule
      height12.3pt depth5.54pt width0pt}%
  \baselineskip=16pt\rm}

\newbox\AuthorBox\newbox\TitleBox
\newbox\TFLinebox
\newbox\FLinebox
\newbox\HLinebox
\def\SetTFLinebox#1{\setbox\TFLinebox=\hbox{#1}}
\def\SetFLinebox#1{\setbox\FLinebox=\hbox{#1}}
\def\SetHLinebox#1{\setbox\HLinebox=\hbox{#1}}

 \def\SetAuthorHead#1{%
     \setbox\AuthorBox=\hbox{\ninepoint \it 
           \ignorespaces\frenchspacing#1\unskip}}
 \def\SetTitleHead#1{%
     \setbox\TitleBox=\hbox{\ninepoint \it
           \ignorespaces\frenchspacing#1\unskip}}

  \def\itSpacing{\relax}
  \def\itSpacingOff{\relax}


 \def\Hrule{\hrule width0pt height0pt}

  \newskip\ProcSkip \ProcSkip 8pt plus2pt minus2pt

 \newskip\LastSkip
 \def\SaveLastSkip{\LastSkip\lastskip}
 \def\RestoreLastSkip{\vskip-\LastSkip\vskip\LastSkip}

 \def\NoindentAfter{\everypar={\setbox0=\lastbox\everypar={}}}

 \long\def\H#1\par#2\par{\notenumber=0 \titlepagetrue%
    {
    \baselineskip=20pt
    \parindent=0pt\parskip=0pt\frenchspacing
    \leftskip=0pt plus .2\hsize minus .3\hsize
    \rightskip=0pt plus .2\hsize minus .3\hsize
 \def\\{\unskip\break}%
    \pretolerance=10000 \Hfont #1\unskip\break
     \vskip7pt\Hrule
\hfill \Authorfont #2\hfill\hfill\unskip}
    \vskip48pt plus 4pt minus 4pt
    \par\NoindentAfter\rm}

 \long\def\Hi#1\par#2\par{\notenumber=0 \titlepagetrue%
    {  \baselineskip=0pt  \parindent=0pt\parskip=0pt\frenchspacing
    \leftskip=0pt plus .2\hsize minus .3\hsize
    \rightskip=0pt plus .2\hsize minus .3\hsize
}
    \rm}


 \newdimen\PageRemainder
  \def\SetPageRemainder{
     \PageRemainder=\pagegoal
     \ifdim\PageRemainder=\maxdimen\PageRemainder=\vsize
     \else\advance\PageRemainder by -1\pagetotal\fi}

  \def\Rpt@{}\def\Rpt@@{}

  \long\def\HH#1\par{\par
  \SaveLastSkip\removelastskip\goodbreak
  \ifdim\LastSkip<30pt 
     \LastSkip 30pt
plus 3pt minus 2pt\fi
  \SetPageRemainder\advance\PageRemainder-\LastSkip
  \ifdim\PageRemainder<150pt
       \edef\Rpt@{remain = \the\PageRemainder\noexpand\\
                pagetotal=\the\pagetotal\noexpand\\
                           pagegoal=\the\pagegoal}%
          \fi
   \ifdim\PageRemainder<65pt 
       \ifdim\PageRemainder > 0pt
          \edef\Rpt@@{\noexpand\\
                      Had HH PageRemainder$<$\relax 65pt\noexpand\\
                      Hence forced break!}%
     \vskip 0pt plus .2\PageRemainder\eject 
    \fi\fi
    \vskip\LastSkip\Hrule 
    \pretolerance=10000\rightskip=0pt plus 3em
    \hangafter1 \hangindent=2.2em%
    \noindent
    \HHfont \unskip \Ednote{\Rpt@\Rpt@@}%
            \def\Rpt@{}\def\Rpt@@{}%
            \ignorespaces
            #1\par\rightskip=0pt\pretolerance=\StdPretolerance%
    \NoindentAfter
\tenpoint\rm%
     \medskip \vskip\ProcSkip}

  \long\def\HHH#1\par{\par%
  \SaveLastSkip\removelastskip\goodbreak
  \ifdim\LastSkip<\ProcSkip%
     \LastSkip\ProcSkip\fi
  \SetPageRemainder\advance\PageRemainder-\LastSkip
  \ifdim\PageRemainder<150pt
       \edef\Rpt@{remain = \the\PageRemainder\noexpand\\
                pagetotal=\the\pagetotal\noexpand\\
                           pagegoal=\the\pagegoal}%
       \fi
   \ifdim\PageRemainder<48pt  
        \ifdim\PageRemainder > 0pt
             \edef\Rpt@@{\noexpand\\
                      Had HHH PageRemainder$<$\relax48pt\noexpand\\
                      Hence forced break!}%
       \vskip 0pt plus .2\PageRemainder\eject 
      \fi\fi
   \vskip\LastSkip\par\noindent
   \HHHfont \unskip\Ednote{\Rpt@\Rpt@@}%
  \def\Rpt@{}\def\Rpt@@{}%
  \ignorespaces
   #1\unskip.\quad\rm\ignorespaces
   \ignorepars}

  \long\def\ignorepars#1\par{\def\Test{#1}%
     \ifx\Test\Empty\def\This{\ignorepars}%
        \else\def\This{\Test\par}\fi
           \This}
  \def\Empty{}

 \def\Abstract#1\par{\bgroup\Smallfonts\narrower\HHH #1\par}
 \def\endAbstract{\par\egroup}


 \def\ProcBreak{\par%
    \ifdim\lastskip<8pt%
    \removelastskip%
    \penalty-200\vskip\ProcSkip\fi}

 \def\th#1\par{\ProcBreak \noindent
   {\thfont\ignorespaces
    #1\unskip.}\it\itSpacing\kern.4em\ignorepars}



 \def\pf#1\par{\ProcBreak %
    \noindent\pffont#1\unskip.\rm\itSpacingOff{\kern .7em}\ignorepars}


  \def\qedbox{\hbox{\vbox{
    \hrule width0.2cm height0.2pt
    \hbox to 0.2cm{\vrule height 0.2cm width 0.2pt
             \hfil\vrule height0.2cm width 0.2pt}
    \hrule width0.2cm height 0.2pt}\kern1pt}}

  \def\qed{\ifmmode\qedbox
    \else\unskip\ \hglue0mm\hfill\qedbox\ProcBreak\fi}

  \def \rk #1\par{\ProcBreak
     \noindent{\rkfont\ignorespaces #1\unskip.}%
     \rm\kern.6em\ignorepars}

  \def \df #1\par{\ProcBreak
     \noindent{\dffont\unskip\ignorespaces #1\unskip.}%
     \rm\kern.6em\ignorepars}

  \def \eg #1\par{\ProcBreak
     \noindent\egfont\unskip\ignorespaces #1\unskip.
     \rm\kern.6em\ignorepars}

  \newdimen\Overhang

   \def\MaxTag@#1#2#3#4#5{\setbox0=\hbox{#4\ignorespaces#2\unskip}%
     \dimen0=\wd0\advance\dimen0 by#3
     \ifdim\dimen0<#5\relax\dimen0=#5\fi
     \expandafter\edef\csname #1Hang\endcsname{\the\dimen0}}

 \def\MaxItemTag#1{\MaxTag@{Item}{#1}{.4em}{\ItemStyle}{\parindent}}%
 \def\MaxItemItemTag#1{%
        \MaxTag@{ItemItem}{#1}{.4em}{\ItemItemStyle}{\parindent}}
 \def\MaxNrTag#1{\MaxTag@{Nr}{#1}{.5em}{\NrStyle}{\parindent}}
 \def\MaxReferenceTag#1{%
        \MaxTag@{Reference}{[#1]}{.6em}{\ninerm}{\parindent}}
 \def\MaxFootTag#1{\MaxTag@{Foot}{#1}{.4em}{\ninerm}{\z@}}

  \def\SetOverhang@{\Overhang=.8\dimen0%
     \advance\Overhang by \wd0\relax
     \ifdim\Overhang>\hangindent\relax
       \advance\Overhang by .25\dimen0%
       \Ednote{Tag is pushing text.}\osumess{Tag is pushing text.}%
     \else\Overhang=\hangindent
     \fi}

   \def\Item#1{\par\noindent
      \hangafter1\hangindent=\ItemHang
      \setbox0=\hbox{\ItemStyle\ignorespaces#1\unskip}%
      \dimen0=.4em\SetOverhang@
      \rlap{\box0}\kern\Overhang\ignorespaces}

   \def\ItemItem#1{\par\noindent
      \hangafter1\hangindent=\ItemItemHang
      \setbox0=\hbox{\ItemItemStyle\ignorespaces#1\unskip}%
      \dimen0=.4em\SetOverhang@
      \advance\hangindent by \ItemHang
      \kern\ItemHang\rlap{\box0}%
      \kern\Overhang\ignorespaces}

  \def\Nr#1{\par\noindent\hangindent=\NrHang 
    \setbox0=\hbox{\NrStyle\ignorespaces#1\unskip}%
    \dimen0=.5em\SetOverhang@
    \rlap{\box0}\kern\Overhang
    \hangindent=\z@\ignorespaces}

   \newskip\Rosterskip\Rosterskip 1pt plus1pt 
   \def\Roster{\par\ifdim\lastskip<\Rosterskip\removelastskip\vskip\Rosterskip\fi
    \bgroup}
   \def\endRoster{\par\global\edef\LastSkip@{\the\lastskip}\removelastskip
       \egroup\penalty-50\LastSkip\LastSkip@\relax
       \ifdim\LastSkip<\Rosterskip\LastSkip\Rosterskip\fi
       \vskip\LastSkip}




 \def\cite#1{
    \def\nextiii@##1,##2\end@{{\frenchspacing\rm 
      \lBr\ignorespaces##1\unskip{\rm,~\ignorespaces##2}\rBr}}%
    \IN@0,@#1@%
    \ifIN@\def\next{\nextiii@#1\end@}\else
    \def\next{{\rm\lBr#1\rBr}}\fi\next}


   \def \Bib#1\par{%
       \par\removelastskip\SetPageRemainder
       \ifdim\PageRemainder < 97pt
        \ifdim\PageRemainder > 0pt
        \vfill\eject
       \fi\fi
    \ProcBreak \par\begingroup\parskip=0 pt%
    \goodbreak \vskip 15 pt plus 10 pt
    \noindent\null\hfill\Bibfont
      \ignorespaces #1\unskip\hfill\null\par 
    \frenchspacing \Smallfonts\rm
    \parskip=2.5 pt plus 1 pt minus.5pt%
    \nobreak\vskip 12pt plus 2pt minus2pt\nobreak
    \leftskip=0 pt \baselineskip=10.5pt}

 \def\ReferenceTagSlide{0em}
  \def\ReferenceTagGap{.5em}

  \def \rf#1{\par\noindent
     \hangafter1\hangindent=\ReferenceHang      
     \setbox0=\hbox{\ninerm[\ignorespaces#1\unskip]}%
     \dimen0=\ReferenceTagGap\SetOverhang@
     \rlap{\kern\ReferenceTagSlide\box0}%
     \kern\Overhang\ignorespaces}

  \def\ref#1\par#2\par#3\par#4\par{%
     \rf{#1}#2\unskip,\ #3\unskip,\
     #4\unskip.}


  \long\def\Coordinates#1\endCoordinates{
 {\par\vskip4pt\def\\{\unskip, }\Coordfont\baselineskip10.5pt\noindent#1}}

 \def\pagecontents{
  \gdef\Pagetot@l{\pagetotal}
  \ifvoid\TRMargIns\else
    \rlap{\kern\hsize\kern10pt\vbox to 0pt{%
         \box\TRMargIns\vss}}\fi
  \ifvoid\topins\else\unvbox\topins\fi
   \dimen@=\dp\@cclv \unvbox\@cclv 
   \ifvoid\footins\else 
     \vskip\skip\footins
     \footnoterule
     \unvbox\footins\fi
   \ifr@ggedbottom \kern-\dimen@ \vfil \fi}


 \newcount\Ht 

 \def \Acc{\expandafter } 

 \def\swthat{\raise -1.1 ex\hbox{\sam$\widehat{}$}}
 \def\swttilde{\raise -1.2 ex\hbox{\sam$\widetilde{}$}}
 \def \overdot{{\raise .2 ex \hbox to 0pt {\hss\bf\smash{.}\hss}}}
 \def \overcircle{{\raise .1 ex \hbox to 0pt
    {\sam$\eightpoint\scriptstyle\hss\circ\hss$}}}

 \def \Mathaccent#1#2{{\sam 
  \setbox4=\hbox{$\vphantom{#2}$}
  \Ht=\ht4 
  \setbox5=\hbox{${#1}$}
  \setbox6=\hbox{${#2}$}
  \setbox7=\hbox to .5\wd6{}
  \copy7\kern .1\Ht \raise\Ht sp\hbox{\copy5}\kern-.1\Ht
  \copy7\llap{\box6}
  }}

  \def\SwtCheck #1{
        \ifmmode \check{#1}%
                \else \v {#1}%
                \fi}

 \def\barpartial {%
   \kern .17 em
    \overline {\kern -.17 em\partial\kern-.03 em}%
    \kern .03 em}

 
  \def\Overline#1{\setbox1=\hbox{\sam ${#1}$}%
      \ifdim \wd1 > 6pt
    \kern .11 em
    \overline {\kern -.11 em#1\kern-.14 em}
    \kern .14 em
  \else
    \kern .03 em
    \overline {\kern -.03 em#1\kern-.04 em}
    \kern .04 em
  \fi}

 \def\SOverline#1{\setbox1=\hbox{\sam ${#1}$}%
      \ifdim \wd1 > 7pt
    \kern .22 em
    \overline {\kern -.22 em#1\kern-.09 em}%
    \kern .09 em
  \else
    \kern .10 em
    \overline {\kern -.10 em#1\kern-.04 em}%
    \kern .04 em
  \fi}


 \def\Underline#1{\setbox1=\hbox{\sam ${#1}$}%
      \ifdim \wd1 > 6pt
    \kern .11 em
    \underline {\kern -.11 em#1\kern-.14 em}
    \kern .14 em
  \else
    \kern .03 em
    \underline {\kern -.03 em#1\kern-.04 em}
    \kern .04 em
  \fi}

 \def\SUnderline#1{\setbox1=\hbox{\sam ${#1}$}%
      \ifdim \wd1 > 7pt
    \kern .04 em
    \underline {\kern -.04 em#1\kern-.2 em}%
    \kern .2 em
  \else
    \kern .0 em
    \underline {\kern -.0 em#1\kern-.15 em}%
    \kern .15 em
  \fi}


 \def \Blackbox
   {\leavevmode\hskip .3pt \vbox
   {\hrule height 5pt\hbox{\hskip 4.5pt}}\hskip .5pt}

 \def \XX{\Blackbox\kern.5pt\Blackbox} 

  \def\.{.\kern1pt}

    \def\Hyphen{\edef\this{\the\hyphenchar\font}%
          \hyphenchar\font=-1\char\this\hyphenchar\font=\this}

 \ifx\undefined\text
  \def\text#1{\hbox{\rm #1}}\fi 



   \everymath{}  

  \def\PassMath@@{\aftergroup\AfterMath@} 

 \let\PassMath@\PassMath@@

 \def\AfterMath@{\futurelet\next\AfterMathMole@}

 \def\AfterMathMole@{
      \ifcat\next\space
          \def\this{}
      \else
      \ifcat\next\egroup %
        \def\this{\osumess{Handset mathsurround?? ---(see dollar brace)}}%
      \else
      \def\this{\AAfterMath@}
      \fi\fi
      \this}

 \def\hyphen@{-}
 \def\paren@{)}
 \def\apostr@{'}

 \def\MSC#1{\kern-.8\mathsurround#1\kern.8\mathsurround}

 \def\AAfterMath@#1{\def\Next{#1}
    \IN@0\Next @,.;:!?\relax @%
    \ifIN@\def\this{\MSC{\Next}}%
    \else
    \ifx\Next\hyphen@\def\this{\futurelet\next\AfterHyphen@}%
    \else
    \ifx\Next\paren@\def\this{#1}%
    \else 
    \ifx\Next\apostr@\def\this{#1}%
    \else \def\this{\osumess{Handset mathsurround??}%
                 #1}\fi\fi\fi\fi
    \this}

 \def\AfterHyphen@#1{\def\Next{#1}%
   \ifx\Next\hyphen@\def\this{--}\else
   \ifcat\next\space%
   \def\this{\kern-\mathsurround\kern.05em- \Next}\else
   \def\this{\kern-\mathsurround\kern.05em\Hyphen\Next}\fi\fi\this}

 \def\sam{\mathsurround=\z@\let\PassMath@\relax}  %
 \def\mas{\mathsurround=\StdMathsurround\let\PassMath@\PassMath@@}
 
 \def\Mas{\mathsurround=\StdMathsurround
                \everymath{\PassMath@}\let\PassMath@\PassMath@@}

 \def\m@th{\mathsurround=\z@\everymath{}}

 \def\m@@th{\mathsurround=\z@\everymath={}\let\m@th\relax}

\def\underbar#1{$\setbox\z@\hbox{#1}\dp\z@\z@
      \m@th \underline{\box\z@}$\relax}

\def\mathhexbox#1#2#3{\leavevmode
  \hbox{\m@@th$\m@th \mathchar"#1#2#3$}}

\def\dots{\relax\ifmmode\ldots\else$\m@th\ldots\,$\relax\fi}

\def\dotfill{\cleaders\hbox{\m@@th$\m@th \mkern1.5mu.\mkern1.5mu$}\hfill}
\def\rightarrowfill{$\m@th\mathord-\mkern-6mu%
  \cleaders\hbox{\m@@th$\mkern-2mu\mathord-\mkern-2mu$}\hfill
  \mkern-6mu\mathord\rightarrow$\relax}
\def\leftarrowfill{$\m@th\mathord\leftarrow\mkern-6mu%
  \cleaders\hbox{\m@@th$\mkern-2mu\mathord-\mkern-2mu$}\hfill
  \mkern-6mu\mathord-$\relax}

\def\downbracefill{$\m@th\braceld\leaders\vrule\hfill\braceru
  \bracelu\leaders\vrule\hfill\bracerd$\relax}
\def\upbracefill{$\m@th\bracelu\leaders\vrule\hfill\bracerd
  \braceld\leaders\vrule\hfill\braceru$\relax}

\def\angle{{\vbox{\m@@th\ialign{$\m@th\scriptstyle##$\crcr
      \not\mathrel{\mkern14mu}\crcr
      \noalign{\nointerlineskip}
      \mkern2.5mu\leaders\hrule height.34pt\hfill\mkern2.5mu\crcr}}}}

\def\big#1{{\m@@th\hbox{$\left#1\vbox to8.5\p@{}\right.\n@space$}}}
\def\Big#1{{\m@@th\hbox{$\left#1\vbox to11.5\p@{}\right.\n@space$}}}
\def\bigg#1{{\m@@th\hbox{$\left#1\vbox to14.5\p@{}\right.\n@space$}}}
\def\Bigg#1{{\m@@th\hbox{$\left#1\vbox to17.5\p@{}\right.\n@space$}}}
\def\n@space{\nulldelimiterspace\z@ \m@th}

\def\root#1\of{\setbox\rootbox\hbox{\m@@th$\m@th\scriptscriptstyle{#1}$}
  \mathpalette\r@@t}
\def\r@@t#1#2{\setbox\z@\hbox{\m@@th$\m@th#1\sqrt{#2}$\relax}
  \dimen@\ht\z@ \advance\dimen@-\dp\z@
  \mkern5mu\raise.6\dimen@\copy\rootbox \mkern-10mu \box\z@}

\def\mathph@nt#1#2{\setbox\z@\hbox{\m@@th$\m@th#1{#2}$}\finph@nt}

\def\mathsm@sh#1#2{\setbox\z@\hbox{\m@@th$\m@th#1{#2}$}\finsm@sh}

\def\@vereq#1#2{\lower.5\p@\vbox{\m@@th\baselineskip\z@skip\lineskip-.5\p@
    \ialign{$\m@th#1\hfil##\hfil$\crcr#2\crcr=\crcr}}}

\def\mathpalette#1#2{\sam\mathchoice{#1\displaystyle{#2}}%
  {#1\textstyle{#2}}{#1\scriptstyle{#2}}{#1\scriptscriptstyle{#2}}\mas}

\def\widehat#1{\setbox\z@\hbox{\sam$#1$}%
 \ifdim\wd\z@>\tw@ em\mathaccent"0\msbfam@5B{#1}%
 \else\mathaccent"0362{#1}\fi}
\def\widetilde#1{\setbox\z@\hbox{\sam$#1$}%
 \ifdim\wd\z@>\tw@ em\mathaccent"0\msbfam@5D{#1}%
 \else\mathaccent"0365{#1}\fi}

 \def\dots{\relax{}
  \ifmmode\def\thedots{\mdots@}\else\def\thedots{\tdots@}\fi %
  \thedots}

 \let\@ldeqno\eqno\let\@ldleqno\leqno
 \def\eqno{\everymath{}\@ldeqno} \def\leqno{\everymath{}\@ldleqno}

  \let\@ldeqalignno\eqalignno
  \def\eqalignno#1{\sam\@ldeqalignno{#1}\mas}
  \let\@ldeqalign\eqalign
  \def\eqalign#1{\sam\@ldeqalign{#1}\mas}

 \def\overrightarrow#1{\vbox{\m@th\ialign{##\crcr
      \rightarrowfill\crcr\noalign{\kern-\p@\nointerlineskip}
      $\hfil\displaystyle{#1}\hfil$\crcr}}}
 \def\overleftarrow#1{\vbox{\m@th\ialign{##\crcr
      \leftarrowfill\crcr\noalign{\kern-\p@\nointerlineskip}
      $\hfil\displaystyle{#1}\hfil$\crcr}}}
 \def\overbrace#1{\mathop{\vbox{\m@th\ialign{##\crcr\noalign{\kern3\p@}
      \downbracefill\crcr\noalign{\kern3\p@\nointerlineskip}
      $\hfil\displaystyle{#1}\hfil$\crcr}}}\limits}
 \def\underbrace#1{\mathop{\vtop{\m@th\ialign{##\crcr
      $\hfil\displaystyle{#1}\hfil$\crcr\noalign{\kern3\p@\nointerlineskip}
      \upbracefill\crcr\noalign{\kern3\p@}}}}\limits}

  \let\@ldmatrix\matrix
  \let\end@ldmatrix\endmatrix
  \def\matrix{\sam\@ldmatrix}
  \def\endmatrix{\end@ldmatrix\mas}
  \let\@ldgather\gather
  \let\end@ldgather\endgather
  \def\gather{\sam\@ldgather}
  \def\endgather{\end@ldgather\mas}
  \let\@ldalign\align
  \let\end@ldalign\endalign
  \def\align{\sam\@ldalign}
  \def\endalign{\end@ldalign\mas}
  \let\@ldaligned\aligned
  \let\end@ldaligned\endaligned
  \def\aligned{\sam\@ldaligned}
  \def\endaligned{\end@ldaligned\mas}
  \let\@ldtag\tag
  \def\tag{\sam\@ldtag}
   %

   \let\MinCDArrowWidth\minCDaw@




\newskip\insertskipamount\newskip\inserthardskipamount
\insertskipamount 6pt plus2pt 
\inserthardskipamount 6pt
\def\insertskip{\vskip\insertskipamount}
\newcount\SplitTest
\def\SetSplitTest{\SplitTest\insertpenalties
  \insert\topins{\floatingpenalty1}%
  \advance\SplitTest-\insertpenalties}
\def\midinsert{\par
 \SaveLastSkip\penalty-150\SetSplitTest\RestoreLastSkip
 \ifnum\SplitTest=-1
  \@midfalse\p@gefalse\else\@midtrue\fi\@ins}
\def\@ins{\par\begingroup\setbox\z@\vbox\bgroup%
  \vglue\inserthardskipamount}
\def\endinsert{\egroup 
  \if@mid \dimen@\ht\z@ \advance\dimen@\dp\z@
    \advance\dimen@\insertskipamount
    \advance\dimen@\pagetotal\advance\dimen@-\pageshrink
    \ifdim\dimen@>\pagegoal\@midfalse\p@gefalse\fi\fi
  \if@mid%
    \ifdim\lastskip<\insertskipamount\removelastskip\insertskip\fi
    \nointerlineskip\box\z@\penalty-200\insertskip
  \else%
    \SaveLastSkip
    \insert\topins{\penalty100 
    \splittopskip\z@skip
    \splitmaxdepth\maxdimen \floatingpenalty\z@
    \ifp@ge \dimen@\dp\z@
    \vbox to\vsize{\unvbox\z@\kern-\dimen@}
    \else \box\z@\nobreak\insertskip\fi}
    \RestoreLastSkip
   \fi\endgroup}


  \newcount\notenumber
  
  \def\note{\advance\notenumber by 1
    \footnote{\the\notenumber)}}

  \newbox\footbox

  \def\footnote#1{\let\@sf\empty
    \ifhmode\edef\@sf{\spacefactor\the\spacefactor}\/\fi
    \sam${}^{\fam0 #1}$\@sf\vfootnote{#1}}%

  \def\vfootnote#1{\insert\footins\bgroup
     \interlinepenalty100 \splittopskip=1pt
     \floatingpenalty=20000
     \leftskip=0pt\rightskip=0pt%
     \parindent=.3em
     \Smallfonts\rm
     \FootItem@{#1}
     \futurelet\next\fo@t}

  \def\FootItem@#1{\par\hangafter1\hangindent=\FootHang
     \setbox0=\hbox{\ignorespaces#1\unskip}%
     \dimen0=.4em\SetOverhang@
     \noindent\rlap{\box0}\kern\Overhang\ignorespaces}


  \def\fo@t{\ifcat\bgroup\noexpand\next \let\next\f@@t
    \else\let\next\f@t\fi \next}
  \def\f@@t{\bgroup\aftergroup\@foot\let\next}
  \def\f@t#1{\baselineskip=10pt\lineskip=1pt
            \lineskiplimit=0pt #1\@foot}%
  \def\@foot{
        \hbox{\vrule height0pt depth5pt width0pt}
        \egroup}
  \skip\footins=12 pt plus 0pt minus 0pt 
  \count\footins=1000 
  \dimen\footins=8in 



 \def\osumess#1{\EdSpider{\immediate\write16{Line \the\inputlineno: #1}}}%
 \def\HideEdStuff{\gdef\EdSpider##1{}}

 \font\BigSym=cmmi10 scaled \magstep 4

 \def\change{\InLMargin{\hbox{\BigSym \char63\kern10pt}}}

 \def\beginchange{\InLMargin{\hbox{\sam\twelvepoint$\heartsuit$\kern10pt}}}

 \def\endchange{\InLMargin{\hbox{\sam\twelvepoint$\spadesuit$\kern10pt}}}

 \def\InLMargin#1{\strut\vadjust{%
     \kern-\strutdepth
     \vtop to \strutdepth{%
         \baselineskip\strutdepth
         \llap{\sam$\smash{\hbox{\EdSpider{#1}}}$}\null}}}

 \def\strutdepth{\dp\strutbox}
 \def\strutheight{\ht\strutbox}

 \def\NoteInRMargin#1{\strut\vadjust{%
     \kern-1.001\strutdepth
     \vtop to \strutdepth{%
       \baselineskip\strutdepth
       \vss\rlap{\ninepoint\unskip\hskip\hsize
         \vtop to 0pt{%
           \hsize=16em\hfuzz=\hsize
           \leftskip=10pt%
           \rightskip=0pt plus 10000pt%
           \baselineskip=9.8pt\lineskip=.2pt%
           \let\\\break
           \noindent\EdSpider{#1}\vss}%
                \kern10pt}\hbox{}}
       }}

 \def\ednote#1{\NoteInRMargin{\tentt #1}}

 \def\cbar{\InLMargin{%
      \dimen0=\strutdepth\advance\dimen0 by \lineskip
      \vrule width 3pt
      height \strutheight depth \dimen0 \kern
      3pt}}

 \def\ccbar{\InLMargin{%
      \dimen0=2\strutdepth\advance\dimen0 by 2\lineskip
      \vrule width 3pt
        height 3\strutheight depth \dimen0 \kern
      3pt}}

 \newinsert\TRMargIns
 \dimen\TRMargIns=\maxdimen

  \def\Ednote#1{\insert\TRMargIns{%
       \vbox to 0pt{\hsize=140pt\hfuzz=\hsize
           \leftskip=6pt%
           \rightskip=0pt plus 10000pt%
           \baselineskip=9.8pt\lineskip=.2pt%
           \let\\\break
           \SetPageRemainder
           \vglue540pt\vglue-\PageRemainder
           \noindent\EdSpider{\tentt #1}\vss}%
       \smallskip}}

 \def\KillEdStuff{\def\ednote##1{}\def\Ednote##1{}%
      \let\change\relax\let\beginchange\relax\let\endchange\relax
       \let\cbar\relax\let\ccbar\relax}


  \topskip=12pt
  \newskip\StdBaselineskip 
  \StdBaselineskip 12pt
  \lineskip=1.1pt
  \lineskiplimit=.8pt
  \widowpenalty=10000 
  \clubpenalty=10000  
  \abovedisplayskip=6pt plus 1pt minus 1pt
  \abovedisplayshortskip=3pt plus 1.5pt
  \belowdisplayskip=6pt plus 1pt minus 1pt
  \belowdisplayshortskip=5pt plus 1pt minus 1pt
  \hfuzz=1.5pt   

  \def\StdPretolerance{100}
  \tolerance=\StdPretolerance

  \newdimen\StdMathsurround
  \StdMathsurround=1.5pt 
  \mathsurround=\StdMathsurround
  \Mas                   

   \def\prose{\relax\hbox{\kern.6\StdMathsurround}}
  
  \def\StdParskip{0pt}    
  \parskip=\StdParskip
  \parindent=0.5cm
 

  \def\Times{ptmr  } 
  \def\TimesI{ptmri  } 
  \def\TimesB{ptmb  }
  \def\TimesBI{ptmbi  }
  \def\HelveticaN{phvrrn }

  =\Times at 10bp
  =\TimesB at 10bp
  \font\tenit=\TimesI at 10bp
  =\TimesBI at 10bp

  \font\tenmrm=cmr10  


    =\Times at 9bp 
    \font\nineit=\TimesI at 9bp 
    =\TimesB at 9bp 
    =\TimesBI at 9bp 

    =\HelveticaN at 9bp 


  =\Times at 12bp
  \font\twelveit=\TimesI at 12bp
  =\TimesB at 12bp


  \font\titleit=\TimesI at 14.4bp
  =\TimesB at 14.4bp

 \SetAuthorHead{AuthorHead} 
 \SetTitleHead{TitleHead}  


  \def\lBr{\raise.125ex\hbox{[\kern.1125ex}}
  \def\rBr{\raise.125ex\hbox{\kern.1125ex]}}

 \setbox\footbox=\hbox{\Smallfonts 2)~}



  \bgroup
  \catcode`\@=11 
  \gdef\itSpacing{%
     \xspaceskip=.31em plus.1em minus.05em \sfcode `f=2001
     \itWarning@\let\itWarning@\itWarning@@}
  \gdef\itSpacingOff{%
     \xspaceskip=0pt \sfcode `f=1000
     \let\itWarning@\relax}
   \global\let\itWarning@\relax
  \gdef\itWarning@@{\errmessage{%
  Special italic spacing already in force
  (you have probably omitted an ``endth'').
  See itSpacing macro in osuPSfnt.sty
         }}
  \egroup

 \fontdimen1\titlebf=0.0pt
 \fontdimen2\titlebf=3.6135pt
 \fontdimen3\titlebf=2.8908pt
 \fontdimen4\titlebf=1.44539pt
 \fontdimen5\titlebf=6.64882pt
 \fontdimen6\titlebf=14.45398pt
 \fontdimen7\titlebf=1.60439pt

 \fontdimen1\tenbi=0.26794pt
 \fontdimen2\tenbi=2.50937pt
 \fontdimen3\tenbi=2.00749pt
 \fontdimen4\tenbi=1.00374pt
 \fontdimen5\tenbi=4.59717pt
 \fontdimen6\tenbi=10.03749pt
 \fontdimen7\tenbi=1.11415pt

 \fontdimen1\twelverm=0.0pt
 \fontdimen2\twelverm=3.01125pt
 \fontdimen3\twelverm=2.409pt
 \fontdimen4\twelverm=1.2045pt
 \fontdimen5\twelverm=5.39615pt
 \fontdimen6\twelverm=12.045pt
 \fontdimen7\twelverm=1.33699pt

 \fontdimen1\twelveit=0.27731pt
 \fontdimen2\twelveit=3.01125pt
 \fontdimen3\twelveit=2.409pt
 \fontdimen4\twelveit=1.2045pt
 \fontdimen5\twelveit=5.37207pt
 \fontdimen6\twelveit=12.045pt
 \fontdimen7\twelveit=1.33699pt

 \fontdimen1\twelvebf=0.0pt
 \fontdimen2\twelvebf=3.01125pt
 \fontdimen3\twelvebf=2.409pt
 \fontdimen4\twelvebf=1.2045pt
 \fontdimen5\twelvebf=5.5407pt
 \fontdimen6\twelvebf=12.045pt
 \fontdimen7\twelvebf=1.33699pt

 \fontdimen1\tenrm=0.0pt
 \fontdimen2\tenrm=2.50937pt
 \fontdimen3\tenrm=2.00749pt
 \fontdimen4\tenrm=1.00374pt
 \fontdimen5\tenrm=4.49678pt
 \fontdimen6\tenrm=10.03749pt
 \fontdimen7\tenrm=1.11415pt

 \fontdimen1\tenit=0.27731pt
 \fontdimen2\tenit=2.50937pt
 \fontdimen3\tenit=2.00749pt
 \fontdimen4\tenit=1.00374pt
 \fontdimen5\tenit=4.47672pt
 \fontdimen6\tenit=10.03749pt
 \fontdimen7\tenit=1.11415pt

 \fontdimen1\tenbf=0.0pt
 \fontdimen2\tenbf=2.50937pt
 \fontdimen3\tenbf=2.00749pt
 \fontdimen4\tenbf=1.00374pt
 \fontdimen5\tenbf=4.61723pt
 \fontdimen6\tenbf=10.03749pt
 \fontdimen7\tenbf=1.11415pt

 \fontdimen1\ninerm=0.0pt
 \fontdimen2\ninerm=2.25842pt
 \fontdimen3\ninerm=1.80673pt
 \fontdimen4\ninerm=0.90337pt
 \fontdimen5\ninerm=4.0471pt
 \fontdimen6\ninerm=9.03374pt
 \fontdimen7\ninerm=1.00273pt

 \fontdimen1\nineit=0.27731pt
 \fontdimen2\nineit=2.25842pt
 \fontdimen3\nineit=1.80673pt
 \fontdimen4\nineit=0.90337pt
 \fontdimen5\nineit=4.02904pt
 \fontdimen6\nineit=9.03374pt
 \fontdimen7\nineit=1.00273pt

 \fontdimen1\ninebf=0.0pt
 \fontdimen2\ninebf=2.25842pt
 \fontdimen3\ninebf=1.80673pt
 \fontdimen4\ninebf=0.90337pt
 \fontdimen5\ninebf=4.15552pt
 \fontdimen6\ninebf=9.03374pt
 \fontdimen7\ninebf=1.00273pt


 \newcount\MaxSpaceFactor
 \MaxSpaceFactor=3000 

 \def\ItemStyle{\rm}
 \def\NrStyle{\rm}
 \def\ItemItemStyle{\rm}

 \MaxItemTag{(iii)}
 \MaxItemItemTag{(iii)}
 \MaxNrTag{(2)}
 \MaxFootTag{2)}
 \def\ReferenceHang{30pt}

 \catcode`\@=\active


\loadbold

=\Times  
=\Times scaled750
=\Times scaled650
\font\rms=\Times scaled 920 

=\TimesBI scaled 860
=\TimesI scaled 860

\textfont0=\rrm  
\scriptfont0=\erm 
\scriptscriptfont0=\srm

\def\Augment#1#2{%
    \toks0\expandafter{#1}\toks2{#2}%
    \edef#1{\the\toks0\the\toks2}}

 \font\twelverma=\Times  scaled 1200
 \font\tenrma=\Times  scaled 1000
 \font\ninerma=\Times scaled 920
 =\Times scaled 840
 \font\sevenrma=\Times scaled 760
 =\Times scaled 680
 \font\fiverma=\Times scaled 600

 \Augment\tenpoint{%
  \textfont0=\tenrma  \scriptfont0=\sevenrma  
  \scriptscriptfont0=\fiverma  }

 \Augment\ninepoint{%
  \textfont0=\ninerma  \scriptfont0=\sevenrma 
  \scriptscriptfont0=\fiverma}

 \Augment\twelvepoint{%
  \textfont0=\twelverma  \scriptfont0=\ninerma  
  \scriptscriptfont0=\sevenrma}

\mathsurround=1pt
\hsize=13.45truecm
\vsize=19.5truecm
\hoffset=1.25truecm
\voffset=2truecm
\advance\baselineskip by 2pt

\predefine\til{\~}
\def\~#1{\relax\ifmmode\widetilde{#1}\else\til{#1}\fi}

\redefine \ge{\geqslant}
\define \wt#1{\mathaccent"0365{#1}}
\define \wh#1{\mathaccent"0362{#1}}

\define\cdvf{complete discrete valuation field }

\define \iss{\,\Mathaccent{\raise -.8 ex\hbox{$\widetilde{}$\kern.1em}}\rightarrow\,}

\define \prlim{{\varprojlim}\vphantom{i}\,}

\define \ur{\mathop{\fam0 ur}}

\define \ab{\mathop{\fam0 ab}}

\define \sep{\mathop{\fam0 sep}}

\define \alg{\mathop{\fam0 alg}}

\define \Frob{\operatorname{\fam0 Frob}}
\define \chr{\mathop{\fam0 char}\,}

\define \Tor{\operatorname{\fam0 Tors}}

\define \Gal{\mathop{\fam0 Gal}}
\define \Hom{\operatorname{\fam0 Hom}}

\Mas
\HideEdStuff
\rm 
 

\def\issn{{\nineit ISSN 1464-8997 (on line) 1464-8989 (printed)}}

\def\gtp{{\nineit Published 10 December 2000: \ \copyright\ Geometry \& 
Topology Publications}}

\def\gtv3{{\nineit Geometry \& Topology Monographs, Volume 3 (2000) --
Invitation to higher local fields}}


\def\lione
{{\rms Geometry \& Topology Monographs}}

\def \litwo{{\rms Volume 3: Invitation to higher local fields
}} 

\def\tinfo #1.#2.#3-#4
{{
\noindent  {\lione} \hfill 
\par 
\vskip-1.5pt
\noindent {\litwo} \hfill
\par 
\vskip-1,5pt
\noindent {\rms Part #1, section #2, pages #3--#4} \hfill
\vskip24pt 
}}

\def\tinfos #1.#2.#3-#4
{{
\noindent  {\lione} \hfill 
\par 
\vskip-1.5pt
\noindent {\litwo} \hfill
\par 
\vskip-1.5pt
\noindent {\rms Pages #3--#4} \hfill
\vskip24pt 
}}

\def\tinfoi #1
{{
\noindent  {\lione} \hfill 
\par 
\vskip-1.5pt
\noindent {\litwo} \hfill
\par 
\vskip-1.5pt
\noindent {\rms Pages iii--xi: Introduction and contents} \hfill
\vskip26pt 
}}


  \def\titlepagehead{\hfil}

  \newif\iftitlepage\titlepagefalse
  \newif\ifblankpage\blankpagefalse
  \def\makeheadline{
     \ifblankpage{}\else%
     \iftitlepage
\vbox{\line{\vbox to 8.5pt{}
\ninerm
\copy\HLinebox \hfill
\hglue5mm\ninebf\folio 
\titlepagehead}}%
      \else
\vbox{\ifodd\pageno\rightheadline\else\leftheadline\fi}%
      \fi\vskip 12pt\fi}%
     \def\rightheadline{\line{\vbox to 8.5pt{}%
      \ninerm
\copy\TitleBox \hfill
\hglue5mm\ninebf\folio}}%
     \def\leftheadline{\line{\vbox to 8.5pt{}%
        \unskip\ninerm\unskip\ninebf\folio\hglue5mm
 \hfill \copy\AuthorBox
}}

 \footline={\ifblankpage{}\else
\iftitlepage\ninepoint\sam\hfill
\line{\vbox to 8.5pt{}
\copy\TFLinebox
\hfill
\hglue5mm 
}
            \else
\ninepoint\sam\hfill
\line{\vbox to 8.5pt{}
\copy\FLinebox
\hfill 
\hglue5mm
}
\hfil\fi\global\titlepagefalse\fi}

\def\blankpage{{\blankpagetrue\noindent\hbox to 10pt{\hss}\vfill
\pagebreak}}

\tenpoint\rm 
 

\SetTFLinebox{\gtp }
\SetFLinebox{\gtv3 }
\SetHLinebox{\issn}

\tinfoi{-3}

\pageno=-3

\Hi{}

{}

\SetAuthorHead{Invitation to higher local fields}
\SetTitleHead{Invitation to higher local fields}

\HH Introduction

\bigskip

This volume is a result of the conference on higher local fields
in M\"unster, August~29--September 5, 1999,
which was supported by SFB 478 ``Geometrische Strukturen in der
Mathematik''.
The conference was organized by I. Fesenko and F. Lorenz.
We gratefully acknowledge great hospitality and
tremendous efforts of Falko Lorenz
which made the conference vibrant.

Class field theory as developed in the first half of
this century 
 is a fruitful generalization and extension of Gauss
reciprocity law; it  describes
abelian extensions of number fields
in terms of objects associated to these fields.
Since its  construction,
one of the important themes of number theory was its generalizations
to other classes of fields or to non-abelian extensions.

In modern number theory one encounters very naturally
schemes of finite type over ${\Bbb Z}$.
A very interesting direction
of  generalization of  class field theory is
to develop a  theory for  higher dimensional   fields ---
finitely generated fields over their prime subfields
(or  schemes of finite type over
$\Bbb Z$ in the geometric language).
Work in this subject, higher (dimensional) class field theory, 
was initiated by A.N. Parshin and K.~Kato independently
about twenty five years ago. For an introduction into several global
aspects of the theory see W. Raskind's review on abelian class field theory
of arithmetic
schemes.

One of the first ideas in higher class field theory is to work with
the Milnor $K$-groups instead of the
multiplicative group in the classical theory.
It is one of the principles of class field theory for number fields
to construct the reciprocity map
 by some blending of class field theories  for local fields.
Somewhat similarly, higher dimensional class field theory
is obtained as a blending of 
 higher dimensional {\it local} class field theories, 
which treat abelian extensions of  {\it higher local fields}.
In this way, the higher local fields were introduced in  mathematics.

A precise definition of higher local fields
will be given in section~1 of Part~I; here
we give  an example.
A complete discrete valuation field $K$ whose residue field is
isomorphic to a usual local field with finite residue field 
is  called a two-dimensional local field.
For example, fields ${\Bbb F}_{p}((T))((S))$,
$\Bbb Q_p((S))$ and
$${\Bbb Q}_{p}\{\!\{T\}\!\}=\left\{ \sum_{-\infty}^{+\infty} a_{i} T^{i} :
a_{i} \in {\Bbb Q}_{p}, \inf v_{p}(a_{i}) > -\infty,
\lim_{i \rightarrow -\infty} v_{p}(a_{i})=+\infty \right\}$$
($v_p$ is the $p$-adic valuation map)
are two-dimensional local fields.
Whereas the first two fields above
can be viewed as  generalizations of functional
 local fields,
the latter field comes in sight as an arithmetical generalization
of $\Bbb Q_p$.

In the classical local case, where $K$ is a complete discrete valuation
field with finite residue field, the Galois group $\Gal(K^{\ab}/K)$
of the maximal abelian
extension of $K$ is approximated by the multiplicative group $K^*$;
and the reciprocity map
$$K^{*} \longrightarrow \Gal(K^{\ab}/K)$$
is close to an isomorphism
(it induces an isomorphism between the  group
 $K^*/N_{L/K}L^*$
and $\Gal(L/K)$ for a finite abelian extension $L/K$,
and it is injective with everywhere dense image).
For two-dimensional local fields $K$ as above, instead of the
multiplicative group $K^{*}$, the Milnor $K$-group
$K_{2}(K)$ (cf.\ Some Conventions
and section~2 of Part~I)
plays an important role. For these fields there is
a reciprocity map
$$K_{2}(K) \longrightarrow \Gal(K^{\ab}/K)$$
which is  approximately an isomorphism
(it induces an isomorphism between  the group $K_2(K)/N_{L/K}K_2(L)$
and $\Gal(L/K)$ for a finite abelian extension $L/K$,
and it has  everywhere dense image; but it is not injective: 
the quotient of $K_2(K)$ by the kernel of the reciprocity map
can be described in terms of topological generators,
see section~6 Part~I).

Similar statements hold in the general case of an $n$-dimensional
local field where one works with the Milnor $K_n$-groups
and their quotients (sections~5,10,11 of Part~I);
and even class field theory of
more general classes of complete discrete valuation
fields can be reasonably developed (sections~13,16 of Part~I).

Since $K_{1}(K)=K^{*}$,  higher local class field theory
contains the classical local class field theory as its one-dimensional
version.

\medskip

The aim of this book is to provide
an introduction to higher local fields  and render the main ideas of this
theory.
The book grew  as an extended version of
 talks given at the conference in M\"unster.
Its expository style aims to introduce the reader into the subject and explain
 main ideas, methods and constructions (sometimes omitting details).
The contributors applied essential efforts to explain
the most important features of their subjects.

Hilbert's words in Zahlbericht that
precious treasures are still hidden in the theory of abelian extensions
are still up-to-date.
We hope that this volume, as the first collection of  main strands
of higher local field theory,  
will be useful
as an introduction and guide on the subject.

\medskip

{\bf The first part} presents the theory of higher local fields,
very often in the more general setting of complete
discrete valuation fields.

\smallskip

Section~1, written by I. Zhukov,  introduces higher local fields 
and 
topologies on their additive and multiplicative groups.
Subsection~1.1 contains all basic definitions
and is referred to in many other sections of the volume.
The topologies  are defined in such a way that the topology of
the residue field is taken into account;
the price one pays is that  multiplication is not
continuous in general, however it is sequentially continuous
which allows one to expand elements into convergent
power series or products.

Section~2,  written by O. Izhboldin, is a short  review of the Milnor $K$-groups
and Galois cohomology groups.
It discusses  $p$-torsion and cotorsion  of the groups
$K_n(F)$ and $K_n^{t}(F)=K_n(F)/\cap_{l\ge1} lK_n(F)$,  an analogue of Satz
90  for the
groups $K_n(F)$ and $K_n^{t}(F)$,
and computation of   $H^{n+1}_m(F)$  where
$F$ is either the rational function field in one variable $F=k(t)$
or the formal power series $F=k((t))$.

Appendix to Section 2, written by M. Kurihara and I. Fesenko, 
contains some basic definitions and properties of differential forms
and Kato's cohomology groups in characteristic $p$
and a sketch of the proof of Bloch--Kato--Gabber's theorem
which   describes
the differential symbol from the Milnor $K$-group $K_n(F)/p$
of a field $F$ of positive characteristic $p$
to the differential module $\Omega_F^n$.

Section~4,  written by J. Nakamura, presents main steps of the proof
of  Bloch--Kato's theorem which states that
the norm residue homomorphism
$$K_q(K)/m\to H^q(K,\Bbb Z/m\, (q))$$
is an isomorphism for a henselian discrete valuation field $K$
of characteristic 0 with residue field of positive characteristic.
This theorem and its proof allows one to simplify 
Kato's original approach to  higher local class field theory.

Section~5,  written by M. Kurihara, is a presentation of  main ingredients of Kato's
higher local class field theory.

Section~6,  written by I. Fesenko, is concerned with certain topologies
on the Milnor $K$-groups of higher local fields $K$
which are related to the topology on the multiplicative group;
their properties are discussed and the structure
of the quotient of the Milnor $K$-groups modulo the intersection
of all neighbourhoods of zero is described.
The latter quotient is called a topological Milnor $K$-group;
it was first introduced by Parshin.

Section~7,  written by I. Fesenko, describes
Parshin's higher local class field theory in characteristic $p$, 
which is relatively easy in comparison with the cohomological approach.

Section~8,  written by S. Vostokov, is a review of known approaches to
explicit formulas
for the (wild) Hilbert symbol not only in the one-dimensional case
but in the higher dimensional case as well.
One of them, Vostokov's explicit formula, is of importance for
the study of topological Milnor $K$-groups
 in section~6 and the existence theorem in section~10.

Section~9,  written by M. Kurihara, introduces his exponential homomorphism
for
a complete discrete valuation field of
characteristic zero, which relates differential forms and
the Milnor $K$-groups of the field, thus helping one
to get an additional information on the structure of the latter.
An application to explicit formulas is discussed in subsection~9.2.

Section~10,  written by I. Fesenko, presents his explicit method to construct
higher local class field theory by using topological $K$-groups
and a generalization of Neukirch--Hazewinkel's axiomatic approaches
to class field theory. Subsection~10.2 presents another simple approach to
class field theory in the 
characteristic $p$ case. 
The case of characteristic 0 is sketched using a concept of
Artin--Schreir trees of extensions (as those extensions
in characteristic 0 which are twinkles of  the characteristic $p$ world).
The existence theorem is discussed in subsection~10.5, being built upon
the results of sections~6 and 8.

Section~11, written by M. Spie\ss , provides a glimpse of
  Koya's and his approach to
the higher local
reciprocity map as a generalization of the classical class formations approach
to the level of complexes of Galois modules.

Section~12,  written by M. Kurihara, sketches his classification of  complete
discrete valuation fields $K$
of characteristic 0 with residue field of characteristic $p$
into two classes depending on the behaviour of the torsion part
of a differential module. For each of these classes, 
subsection~12.1 characterizes the quotient filtration of the Milnor
$K$-groups of $K$, for all sufficiently large members of the filtration, 
as  a quotient of differential modules.
For a higher local field the previous result and
higher local class field theory imply certain restrictions
on types of cyclic extensions of the field of sufficiently large degree.
This is described in~12.2.

Section~13,  written by M. Kurihara, describes his  theory of cyclic
$p$-extensions of
an absolutely unramified complete discrete valuation field $K$
with {\it arbitrary} residue field of characteristic $p$.
In this theory a homomorphism is constructed from
 the $p$-part of the group of characters of $K$
to Witt vectors over its residue field.
This homomorphism satisfies some important properties listed
in the section. 

Section~14,  written by I. Zhukov, presents some explicit methods of
constructing abelian extensions
of complete discrete valuation fields.
His approach to  explicit
equations of a cyclic extension of degree $p^n$
which contains a given cyclic extension of degree $p$
is explained.
An application to the structure of topological $K$-groups
of an absolutely unramified higher local field is given in subsection~14.6.

Section~15,  written by J. Nakamura, contains a list of all known results
 on the quotient filtration on the Milnor $K$-groups (in terms of
differential
forms of the residue field) of
a complete discrete valuation field. It discusses his recent study of
the case of a tamely ramified field of characteristic 0 with residue field
of
characteristic $p$ by
using  the exponential map of section~9 and a syntomic complex.

Section~16,  written by I. Fesenko, is devoted to his generalization of
one-dimensional class field theory to a description of
abelian totally ramified  $p$-extensions of
a complete discrete valuation field with arbitrary non separably-$p$-closed
residue field. 
In particular, subsection~16.3 shows that two such extensions
coincide if and only if
their norm groups coincide.
An illustration to the theory of section~13 is given in subsection~16.4.

Section~17,  written by I. Zhukov, is a review of his recent approach to
ramification theory of a complete discrete valuation field
with residue field whose $p$-basis consists of at most one element.
One of important ingredients of the theory is Epp's theorem
on elimination of wild ramification (subsection~17.1).
New lower and upper filtrations are defined
(so that cyclic extensions of degree $p$ may have
non-integer ramification breaks, see examples in subsection~17.2).
One of  the advantages of this theory is its compatibility
with the reciprocity map.  
A refinement of the filtration for two-dimensional local fields
which is compatible with the reciprocity map is discussed.

Section~18,  written by L. Spriano, presents ramification theory of monogenic
extensions
of complete discrete valuation fields;
his recent study demonstrates that in this case
there is a satisfactory theory if one
systematically uses a generalization of the function $i$ and not $s$
(see subsection~18.0 for  definitions).
Relations to Kato's conductor are discussed in~18.2 and 18.3.

These  sections 17 and 18 can be viewed as 
the rudiments of higher ramification theory; 
 there are several other  approaches. 
Still, there is  no  satisfactory general
ramification theory  for complete discrete valuation fields
in the imperfect residue field case;
 to construct such a theory is a challenging problem.

\smallskip

Without attempting to list all links between the sections
we just mention several paths (2 means Section 2 and Appendix to Section 2)
$$
\alignat2
&1\to 6 \to 7\qquad &&\text{\rm (leading to Parshin's approach in positive
characteristic)},\\
&2 \to 4 \to 5\to 11\qquad &&\text{\rm
(leading to Kato's cohomological description}\\
&\phantom{2 \to 4 \to 5\to 11}\qquad&&\text{\rm of the reciprocity map
and generalized class formations)},\\
&8.3\to 6\to 10\qquad&&\text{\rm (explicit construction of the reciprocity
map)},\\
& 5\to 12\to 13\to 15, \qquad &&\text{\rm (structure of the Milnor
$K$-groups of the fields}\\
&  1\to 10\to 14, 16 \qquad&&\text{\rm and
more explicit study of abelian extensions)},\\
&8,9\qquad &&\text{\rm (explicit formulas for the Hilbert norm symbol}\\
& &&\text{\rm and its
generalizations)},\\
&1\to 10\to 17, 18\qquad&&\text{\rm (aspects of higher ramification
theory)}.
\endalignat
$$

\medskip

A special place in this volume (between Part I and Part II) is occupied
by the work of K.~Kato 
on the existence theorem in higher local class field theory
which  was produced in 1980 as an IHES preprint and has never been published.
We are grateful to K.~Kato for his permission to include
this work in the volume. 
In it, viewing higher local fields as ring objects in the category of iterated
pro-ind-objects, a definition of open subgroups in the Milnor $K$-groups
of  the fields is given.
The self-duality of the additive group of a higher local field is proved.
By studying norm groups of cohomological objects and using cohomological
approach
to higher local class field theory the existence theorem is proved. 
An alternative approach to 
the description of norm subgroups of Galois extensions of higher local
fields and the existence theorem is contained in sections~6 and~10.

\medskip

{\bf The second part} is concerned with various applications and connections
of higher local fields with several other areas.

\smallskip

Section~1, written by A.N. Parshin, describes some first steps in
extending Tate--Iwasawa's analytic method
to define an $L$-function
in higher dimensions; historically the latter problem
 was one of the stimuli of the work  on higher class field theory. 
For generalizing this method
the author advocates the usefulness  of 
the classical Riemann--Hecke approach (subsection~1.1), 
his adelic complexes (subsection~1.2.2) together with his generalization
of Krichever's correspondence (subsection~1.2.1).
He analyzes dimension~1 types of functions in subsection~1.3
and discusses properties of the lattice of commensurable classes
of subspaces in the adelic space associated to a divisor on an algebraic surface
in subsection~1.4. 

Section~2, written by D. Osipov, is a review of his recent
work on adelic constructions of direct images of differentials and symbols
in the two-dimensional case in the relative situation.
In particular, reciprocity laws for relative residues of differentials and
 symbols are introduced and applied to a construction
of the Gysin map for Chow groups.

Section~3, written by A.N. Parshin, presents his theory
of Bruhat--Tits buildings over higher dimensional local fields.
The theory is illustrated with the buildings
for $PGL(2)$ and $PGL(3)$ for one- and two-dimensional local fields.

Section~4, written by E.-U. Gekeler,
provides a survey of  relations between Drinfeld modules
and higher dimensional fields of positive characteristic.

Section~5, written by M. Kapranov,
sketches his recent approach to elements of harmonic analysis on algebraic groups
over functional two-dimensional local fields.
 For a two-dimensional local field
subsection~5.4 introduces 
a  Hecke algebra
which  is formed by operators
which integrate pro-locally-constant complex functions
over a non-compact domain.

Section~6, written by L. Herr, 
is a survey of his recent study of applications of   Fontaine's theory
of
$p$-adic representations of local fields ($\Phi-\Gamma$-modules)
to Galois cohomology of local fields
and explicit formulas for the Hilbert symbol (subsections~6.4--6.6).
The two Greek letters lead to two-dimensional local objects (like 
$\Cal O_{\Cal E(K)}$ introduced in subsection~6.3).

Section~7, written by I. Efrat, 
introduces recent advances in the zero-dimensional anabelian
geometry,  
that is a characterization of fields by means of their absolute Galois
group
(for finitely generated fields and for higher local fields).
His method of construction of henselian valuations on fields
which satisfy some $K$-theoretical properties
is presented in subsection~10.3, 
and applications to an algebraic proof of the local correspondence
part of Pop's theorem and to
higher local fields are given.

Section~8, written by A. Zheglov, presents his  
study of two dimensional local skew fields which was initiated by
A.N. Parshin.
If the skew field has one-dimensional residue field
which is in its centre, 
then  one is naturally led to the study of
automorphisms of the residue field which are associated to a local parameter
of
the skew field.
Results on such automorphisms are described in subsections~8.2 and 8.3.

Section~9, written by I. Fesenko, is an exposition 
of his recent work on noncommutative local reciprocity maps
for totally ramified Galois extensions
with arithmetically profinite group
(for instance $p$-adic Lie extensions).
These maps in general are not homomorphisms but  Galois cycles;
a description of their image and kernel is included.

Section~10, written by B. Erez, is a concise
survey of Galois module theory links with class field theory;
 it lists several open problems.

\medskip

The theory of higher local fields has several interesting aspects
and applications which are not contained in this volume.
One of them is the work of Kato on applications of
an explicit formula for the reciprocity map 
in higher local fields
to calculations of  special values of the $L$-function of a
modular form.
There is some interest in two-dimensional local fields
(especially of the functional type) in certain parts of
mathematical physics, infinite group theory and topology 
where formal power series objects play a central role.

\medskip

Prerequisites for most sections in the first part of the book
are  small:
local fields and local class field theory,
for instance, as presented in Serre's ``Local Fields'',
Iwasawa's ``Local Class Field Theory''
or Fesenko--Vostokov's ``Local Fields and Their Extensions''
(the first source contains a cohomological approach
whereas the last two are cohomology free)
and some basic knowledge of  Milnor $K$-theory
of discrete valuation fields (for instance Chapter~IX of
the latter book).
See also Some Conventions and Appendix to Section~2 of Part~I where we
explain several
notions useful for reading Part~I.

\medskip

We thank P. Schneider for his support of the conference and work on
 this volume. The volume is typed using a modified version of osudeG
 style (written by Walter Neumann and Larry Siebenmann and available
 from the public domain of Department of Mathematics of Ohio State
 University, pub/osutex); thanks are due to Larry for his advice on
 aspects of this style and to both Walter and Larry for permission to
 use it.

\bigskip
Ivan Fesenko \quad Masato Kurihara
\hfill September 2000

\vfill\eject

\pageno=-11

\HH Contents

\quad Some Conventions \dotfill  1

\smallskip

\noindent \quad Part I \dotfill 3

\smallskip

\noindent 1. Higher dimensional local fields (I. Zhukov)  \dotfill  5

\noindent 2. $p$-primary part of the Milnor $K$-groups and Galois cohomology of fields of 

\noindent characteristic $p$  
(O. Izhboldin)  \dotfill   19

\noindent A. Appendix to Section 2 
(M. Kurihara and I. Fesenko) \dotfill 31

\noindent 4. Cohomological symbol  for henselian discrete valuation fields 
of 

\noindent mixed characteristic 
(J. Nakamura)  \dotfill  43 

\noindent 5. Kato's higher local class field theory (M. Kurihara)  
\dotfill  53  

\noindent 6. Topological Milnor $K$-groups of higher local fields (I. Fesenko)  \dotfill   61

\noindent 7. Parshin's higher local class field theory in characteristic $p$
(I. Fesenko)  \dotfill  75  

\noindent 8. Explicit formulas for the Hilbert symbol (S. V. Vostokov)  \dotfill  81

\noindent 9. Exponential maps and explicit formulas (M. Kurihara)  \dotfill   91

\noindent 10. Explicit higher local class field theory (I. Fesenko)  \dotfill  95  

\noindent 11. Generalized class formations and higher class field theory
(M. Spie\ss)  \dotfill   103

\noindent 12. Two types of complete discrete valuation fields
(M. Kurihara)   \dotfill  109

\noindent 13. Abelian extensions of absolutely unramified complete discrete   

\noindent valuation fields (M. Kurihara)  \dotfill  113 

\noindent 14. Explicit abelian extensions of complete discrete valuation fields
(I. Zhukov)  \dotfill  117 

\noindent 15. On the structure of the Milnor $K$-groups of complete discrete valuation 

\noindent fields (J. Nakamura)  \dotfill  123 

\noindent 16. Higher class field theory without using $K$-groups
(I. Fesenko)  \dotfill   137

\noindent 17. An approach to higher ramification theory (I. Zhukov)  \dotfill  143

\noindent 18. On ramification theory of monogenic extensions (L. Spriano) \dotfill  151

\smallskip

\noindent \quad Existence theorem for higher local class field theory (K.~Kato) \dotfill 
165

\smallskip

\noindent \quad Part II \dotfill 197

\smallskip

\noindent 1. Higher dimensional local fields and $L$-functions 
(A. N. Parshin) \dotfill  199

\noindent 2. Adelic constructions for direct images of differentials and
symbols (D. Osipov) \dotfill  215 

\noindent 3. The Bruhat--Tits buildings over 
higher dimensional local fields (A. N. Parshin) \dotfill  223

\noindent 4.  Drinfeld modules and 
local fields of positive characteristic
(E.-U. Gekeler) \dotfill   239

\noindent 5. Harmonic analysis on algebraic groups over
two-dimensional local fields of 

\noindent equal  
characteristic 
(M. Kapranov) \dotfill 
255

\noindent 6. $\Phi$-$\Gamma$-modules and Galois cohomology
(L. Herr) \dotfill   263

\noindent 7. Recovering higher global and local fields 
from Galois groups --- 

\noindent an algebraic approach 
(I. Efrat) \dotfill   273

\noindent 8. Higher local skew fields (A. Zheglov)
\dotfill   281

\noindent 9. Local reciprocity cycles (I. Fesenko) 
\dotfill   293

\noindent 10. Galois modules and class field theory 
(B. Erez) \dotfill 299

\vfill\eject

\pageno=1

\SetAuthorHead{Invitation to higher local fields}
\SetTitleHead{Invitation to higher local fields}

\SetTFLinebox{\gtp }
\SetFLinebox{\gtv3 }
\SetHLinebox{}

\HH Some Conventions

\phantom{}\par

\phantom{}\par

\par 

The notation $X\subset Y$ means that $X$ is a subset of $Y$.

For an abelian group $A$ written additively
denote by $A/m$ the quotient group $A/mA$ where $mA=\{ma:a\in A\}$
and by ${}_mA$ the subgroup of elements of order dividing $m$.
The subgroup of torsion elements of $A$ is denoted by $\Tor A$. 

For an algebraic closure $F^{\alg}$ of $F$ 
denote the separable closure of the field $F$ 
by $F^{\sep}$; let $G_F=\Gal(F^{\sep}/F)$ be the absolute Galois group
of $F$.  
Often for a $G_F$-module $M$
we write
$H^i(F,M)$ instead of $H^i(G_F,M)$.

For a positive integer $l$ which is prime to characteristic of $F$
(if the latter is non-zero) 
denote by $\mu_l=\langle\zeta_l\rangle$ the group of $l$th roots of unity in  $F^{\sep}$.

If $l$ is prime to $\chr(F)$, for $m\ge0$ denote by
$\Bbb Z/l(m)$ the $G_F$-module 
$\mu_{l}^{\otimes m}$ and put $\Bbb Z_l(m)=\prlim_r \, \Bbb Z/l^r(m)$; 
for $m<0$ put 
$\Bbb Z_l(m)=\Hom(\Bbb Z_l, \Bbb Z_l(-m))$.

Let $A$ be a commutative ring.
The group of invertible elements of $A$ is denoted by $A^*$. 
Let $B$ be an $A$-algebra.
$\Omega_{B/A}^1$ denotes as usual the $B$-module of regular differential
forms of $B$ over $A$; $\Omega_{B/A}^n=\wedge^n \Omega_{B/A}^1$.
In particular, $\Omega_A^n=\Omega_{A/\Bbb Z 1_A}^n$
where $1_A$ is the identity element of $A$ with respect
to  multiplication. 
For more on differential modules see subsection A1 of
the appendix to the section~2 in the first part.  


Let $K_{n}(k)=K_{n}^{M}(k)$ be the Milnor $K$-group of a field $k$
(for the definition see subsection~2.0 in the first part). 

For a  \cdvf  $K$ denote by 
$\Cal O= \Cal O _{K}$ its {\it ring of integers},
by $\Cal M=\Cal M_K$ the {\it maximal ideal} of $\Cal O$ 
and by $k=k_{K}$ its {\it residue
field}.
If $k$ is of characteristic $p$,
 denote by $\Cal R$ the set of
{\it Teichm\"{u}ller representatives} (or {\it multiplicative representatives})
in $\Cal O$. 
For $\theta$ in the maximal perfect subfield of $k$ 
 denote by $[\theta]$  its
Teichm\"{u}ller representative.

For a field $k$ denote by $W(k)$  the ring of Witt vectors
(more precisely, Witt $p$-vectors
where $p$ is a prime number)  over $k$.
Denote by $W_r(k)$ the ring of Witt vectors of length $r$ over $k$. 
If $\chr(k)=p$ 
denote by ${\bold F}\colon W(k)\to W(k)$, 
${\bold F}\colon W_r(k)\to W_r(k)$ the map 
$(a_0,\dots)\mapsto (a_0^p,\dots)$. 

Denote by $v_K$ the surjective discrete valuation
$K^*\to\Bbb Z$ (it is sometimes called the {\it normalized discrete valuation} 
of $K$). 
Usually $\pi=\pi_K$ denotes a {\it prime element} of $K$:
$v_K(\pi_K)=1$.

Denote by $K_{\ur}$ the {\it maximal unramified extension} of $K$.
If $k_K$ is finite, denote by $\Frob_K$ the {\it Frobenius automorphism} of
$K_{\ur}/K$.

For  a finite extension $L$ of a complete discrete valuation field
$K$  $\Cal D_{L/K}$  denotes its different.

If $\chr(K)=0$, $\chr(k_K)=p$, then $K$ is called
a field of {\it mixed characteristic}.
If $\chr(K)=0=\chr(k_K)$, then $K$ is called a field of {\it equal characteristic}.

If $k_K$ is perfect, $K$ is called a {\it local field}.

\vfill
\pagebreak 
\end